\documentclass[%
pof,
numerical, 
amsmath,
amssymb,
reprint,
]{revtex4-1}
\usepackage{hyperref}
\usepackage{booktabs}
\usepackage{graphicx}
\usepackage{subfigure}
\usepackage{float}
\usepackage[export]{adjustbox}
\usepackage{dcolumn}
\usepackage{bm}
\usepackage{cleveref}
\usepackage[abbreviations]{glossaries-extra}
\usepackage[utf8]{inputenc}
\usepackage[T1]{fontenc}
\usepackage{threeparttable}
\usepackage{mathptmx}
\usepackage{siunitx}
\usepackage{etoolbox}
\usepackage{color}
\usepackage{CJKutf8}
\usepackage{xcolor}
\usepackage{comment}
\usepackage{multirow}
\usepackage{soul}

\soulregister\cite7
\soulregister\ref7
\soulregister\eqref7
\soulregister\pageref7 

\usepackage{todonotes}

\makeatletter
\def\@email#1#2{%
 \endgroup
 \patchcmd{\titleblock@produce}
  {\frontmatter@RRAPformat}
{\frontmatter@RRAPformat{\produce@RRAP{*#1\href{mailto:#2}{#2}}}\frontmatter@RRAPformat}
  {}{}
}%
\makeatother
\begin{document}

\preprint{AIP/123-QED}

%
\title[]{Linear stability of lattice Boltzmann models with non-ideal equation of state}
\author{S.A.~Hosseini}
\affiliation{Computational Kinetics Group, Department of Mechanical and Process Engineering, ETH Zürich, 8092 Zurich, Switzerland}%
\author{I.V.~Karlin}
\affiliation{Computational Kinetics Group, Department of Mechanical and Process Engineering, ETH Zürich, 8092 Zurich, Switzerland}%
\email{shosseini@ethz.ch, ikarlin@ethz.ch}
\date{\today}

\begin{abstract}
{Detailed study of spectral properties and of linear stability is presented for a class of lattice Boltzmann models with a non-ideal equation of state.
Examples include the van der Waals and the shallow water models.
Both analytical and numerical approaches demonstrate that linear stability requires boundedness of propagation speeds of normal eigen-modes.
The study provides a basis for the construction of unconditionally stable lattice Boltzmann models.
}
\end{abstract}

\maketitle
\section{Introduction}
\label{sec:Introduction}
A variety of flows rely on nonideal equations of state featuring a nonlinear dependence of the thermodynamic pressure on fluid density. These include, but are not limited to, the cubic class of equations of state used to model non-ideal fluids with multiple phases \cite{hosseini2023lattice} and quadratic equations of state featured by shallow water equations \cite{tan1992shallow}. The nonlinearity in the equations of state is directly translated into the propagation speed of the normal eigen-modes commonly referred to as the speed of sound. This in turn comes with challenges to be considered when developing numerical models for the corresponding continuum-level balance equations.

Here, we are interested in a numerical method introduced in the late 80's to solve the Navier--Stokes (NS) equations in the incompressible limit: The lattice Boltzmann method (LBM). The LBM has been successful in modeling flows in the incompressible limit, most notably due to its efficiency, brought about by the low intensity of numerical operators and locality of closure for the pressure, as compared to alternatives relying on pressure closures such as the Poisson equation. Recent years have also witnessed successful extension of the LBM to modeling of fluids with non-ideal equations of states \cite{chen2014critical,hosseini2023lattice}. The so-called {free energy} and {pseudo-potential} LBM are good examples of this growing success \cite{shan1993lattice,swift1996lattice,hosseini2022towards}. However, as with the original LBM for ideal fluids, stability is an important topic of discussion \cite{Geier_Cumulant,malaspinas_recursivereg,ansumali2003minimal,ansumali2002single,Ali_reviewEntropic}. Contrary to the LBM with ideal equation of state where there is an extensive literature on stability analysis and spectral properties, non-ideal LBM literature is devoid of such studies.

In this work, we present for the first time a comprehensive linear analysis of spectral properties and stability of the LBM with a non-ideal equation of state. Being interest in the design of stable attractors guaranteeing unconditional linear stability, as opposed to control over relaxation path, we only consider the lattice Bhatnagar--Gross--Krook (LBGK) model \cite{bhatnagar1954model}. 

{The outline of the paper is as follows: A brief introduction to the LBGK formalism and a discussion of non-linear equations of state are provided in Sec.\ \ref{sec:LBGK}.
In Sec.\ \ref{sec:D1Q3}, we report a comprehensive analysis of linear stability of  the one-dimensional LBGK system. This analysis reveals the necessary and sufficient stability conditions as a bound on the normal eigen-modes. These results extend our previous findings \cite{PhysRevE.110.015306} pointing at the need for asymptotic freedom for unconditional linear stability.
The analysis is illustrated by two typical equations of state, the van der Waals and the shallow water systems.
The one-dimensional LBM analysis is extended to higher dimensions in Sec.\ \ref{sec:D2Q9}. Both analytical and numerical studies of spectral dispersion and dissipation of hydrodynamic modes, and evaluation of the linear stability domain are performed for the two-dimensional LBGK. In particular, we demonstrate that the linear stability condition of normal eigen-modes  remains necessary also in higher dimensions.
The paper closes in Sec.\ \ref{sec:conslusions} with a discussion on the conditions for the construction of unconditionally linear-stable equilibria by introduction of velocity-dependent  pressure in the spirit of the entropic equilibrium for the isothermal ideal equation of state \cite{PhysRevE.110.015306,ansumali2003minimal,ansumali2002single}. }
\section{Generic LBGK \label{sec:LBGK}}
We consider the LBGK model~\cite{qian1992lattice} for nearly-incompressible flows,
\begin{equation}
	f_i(\bm{r}+\bm{c}_i \delta t, t+\delta t) - f_i(\bm{r}, t)= 2\beta\left( f_i^{\rm eq}(\rho, \bm{u}) - f_i\right).\label{eq:LBGK}
\end{equation}
Here $f_i$ are the populations of the discrete velocities $\bm{c}_i$, $i=1,\dots, Q$, $\bm{r}$ is the position in space, $t$ is the time, $\delta t$ is the time-step, $\rho$ is the fluid density and $\bm{u}$ is the flow velocity,
\begin{align}
	&\sum_{i=1}^Q \{1,\bm{c}_i\}f_i = \{\rho,\rho\bm{u}\}.
\end{align}
We consider the standard first-neighbor lattices in space dimension $D$ defined as a $D$-fold tensor product of one-dimensional velocities $c_{i\alpha}\in\{-1,0,1\}$.
For these D$D$Q$3^D$ lattices, a generic class of equilibria $f_i^{\rm eq}$ is defined as follows: First, we introduce a triplet of functions $\Psi_{i\alpha}(u,\mathcal{P})$, $i\alpha\in\{-1,0,1\}$: 
\begin{align}
&\Psi_0=1-\mathcal{P},\\
&\Psi_{-1}=\frac{1}{2}(-u+\mathcal{P}),\\
&\Psi_{1}=\frac{1}{2}(u+\mathcal{P}). 
\end{align}
For a $D$-dimensional lattice, equilibria are defined by a product form,
    \begin{equation}\label{eq:generic_isothermal_EDF}
	f_i^{\rm eq} (\rho,\bm{u})=  \rho  \prod_{\alpha=1}^{D} \Psi_{i\alpha}(u_\alpha,\mathcal{P}^{\rm eq}_{\alpha\alpha}).
\end{equation}
Here 
$\mathcal{P}^{\rm eq}_{\alpha\alpha}$ are diagonal components of the equilibrium pressure tensor divided by the density.
\begin{align}
	\mathcal{P}^{\rm eq}_{\alpha\alpha}= {\pi^*(\rho)} + u_\alpha^2.\label{eq:Peq_gen}
\end{align}
Furthermore, $\beta\in[0,1]$ is the relaxation parameter which is tied to the kinematic viscosity,
\begin{align}
\nu=\pi^* \delta t\left(\frac{1}{2\beta}-\frac{1}{2}\right),\label{eq:viscosity}
\end{align}
and $\varsigma=c_s\delta r/\delta t$ is the lattice speed of sound, $\delta r$ is the lattice spacing while $c_s$ is a pure constant dependent on the choice of the lattice. Below, we use lattice units by setting $\delta r=\delta t=1$. The D$D$Q$3^D$ lattices are characterized by the lattice speed of sound, 
\begin{equation}
		\varsigma = \frac{1}{\sqrt{3}}.\label{eq:cs_lattice}
\end{equation}
The LBGK setup becomes complete once the function $\pi^*$ is specified. Note that, at variance with our previous study~\cite{PhysRevE.110.015306} where only the ideal gas equation of state was considered, in this paper we address a generic non-linear equation of state,
\begin{equation}
    \partial_\rho \pi^* \neq 0.
\end{equation}
The absence of a forcing term in the LBGK equation \eqref{eq:LBGK} points to that we only consider models where the pressure is introduced via the equilibrium function.

Throughout the paper, in order to illustrate general statements, we will make use of two non-ideal equations of state:\\ 
\noindent (a) The shallow water equation of state defined as,
\begin{equation}\label{eq:swe}
    \pi^* = \frac{g}{2}\rho,
\end{equation}
where $g$ is the gravitational constant. Note that, in the context of the shallow water model, the field $\rho$ refers to the water column elevation.

\noindent (b) The second equation of state used for illustration purposes in this paper is that of van der Waals,
\begin{equation}\label{eq:vdw}
    \pi^* = \frac{R T}{1 - b\rho} - a\rho,
\end{equation}
where $R$ is the gas constant, $T$ the temperature and $a$ and $b$ the long-range attractive force and excluded volume parameters. For the van der Waals equation of state, different from shallow water, not all temperature/density combinations lead to a hyperbolic system with finite speed eigen-modes.
\begin{figure}[h]
	\centering
		\includegraphics[width=0.8\linewidth]{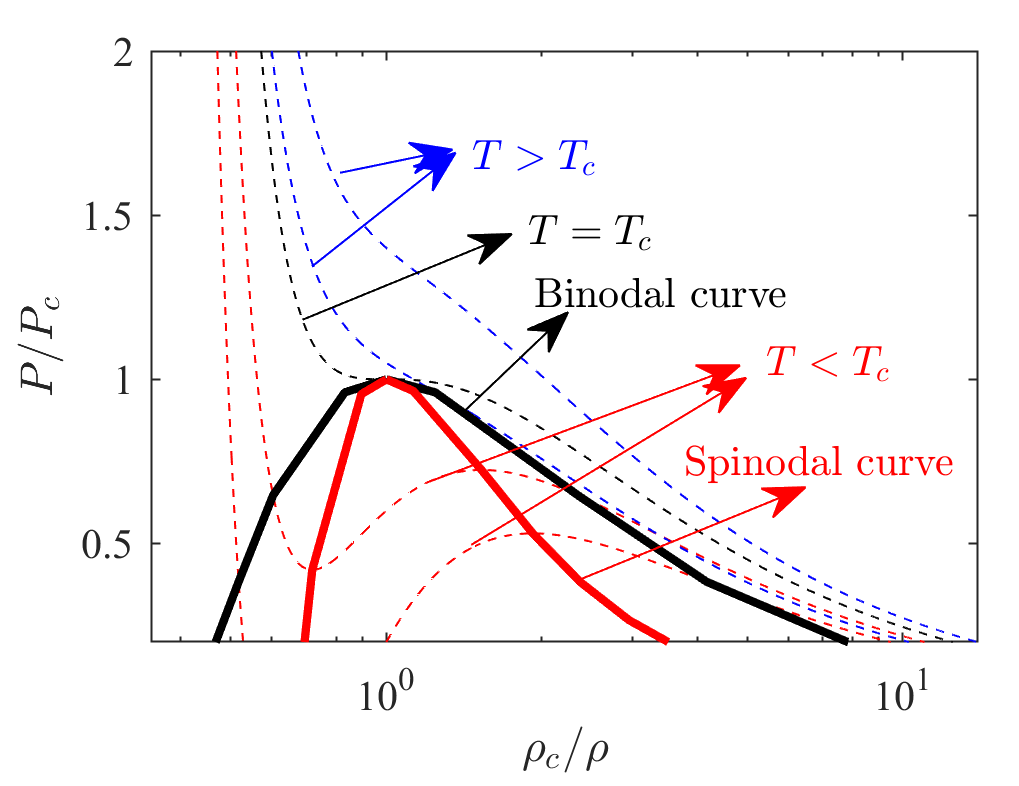}
	\caption{Clapeyron diagram of van der Waals equation of state. Here $\rho_c$ and $P_c$ are the critical density and pressure.}
	\label{clapyeron}
\end{figure}
A typical van der Waals $(\rho\pi^*,1/\rho)$ diagram is illustrated in Fig.~\ref{clapyeron}. Moving on the isobars shown in Fig.~\ref{clapyeron} it is clear that the positive slopes within the spinodal area would lead to negative values of $\partial_\rho (\rho\pi^*)$ which would result in a non-hyperbolic system. Below, the van der Waals equation will make use of values of density on the saturation curve, either the vapor or the liquid branch. The vapor-liquid co-existence densities for a given temperature are computed using the Maxwell equal area construction, see e.\ g. \cite{hosseini2023lattice} for more details.
\section{Linear stability: The D$1$Q$3$ system\label{sec:D1Q3}}
\subsection{Necessary stability condition: Small lattice Knudsen number}
As a first step we consider the one-dimensional LBGK system on the D$1$Q$3$ lattice. In this case, we shall be able to derive in-full the linear stability conditions analytically.
We follow the analytical approach proposed in \cite{Ali_StabilityEntropic,PhysRevE.110.015306}, focusing first on stability in the hydrodynamic limit, which is a necessary condition for the stability of the full discrete system.

At the Euler level, the system is described via the following balance equations for mass and momentum,
\begin{align}
    &\partial_t \rho + \partial_x (\rho u) = 0,\label{eq:continuityE}\\
    &\partial_t (\rho u) + \partial_x \left(\rho u^2+\rho\pi^*\right) = 0.\label{eq:momentumE}
\end{align}
Eigen-mode analysis of the Euler system \eqref{eq:continuityE} and \eqref{eq:momentumE}, reveals the propagation speeds of the two eigen-modes,
\begin{align}\label{eq:hydro_modes_D1Q3}
   { c^{\pm}= u 	\pm \varsigma_\rho,}
\end{align}
{where $\varsigma_\rho$ is the speed of sound,
\begin{equation}\label{eq:sound}
\varsigma_\rho=\sqrt{\partial_\rho(\rho\pi^*)}.
\end{equation}
Speed of sound \eqref{eq:sound} is well-defined, provided 
\begin{equation}
    \partial_\rho(\rho\pi^*)\ge 0,\label{eq:sound_existence}
\end{equation}
which we always assume below. Moreover, we use the convention of positive square root in \eqref{eq:sound} which leads to $c^+\geq c^-$.}

\begin{figure}[h!]
	\centering
	\includegraphics[width=0.7\linewidth]{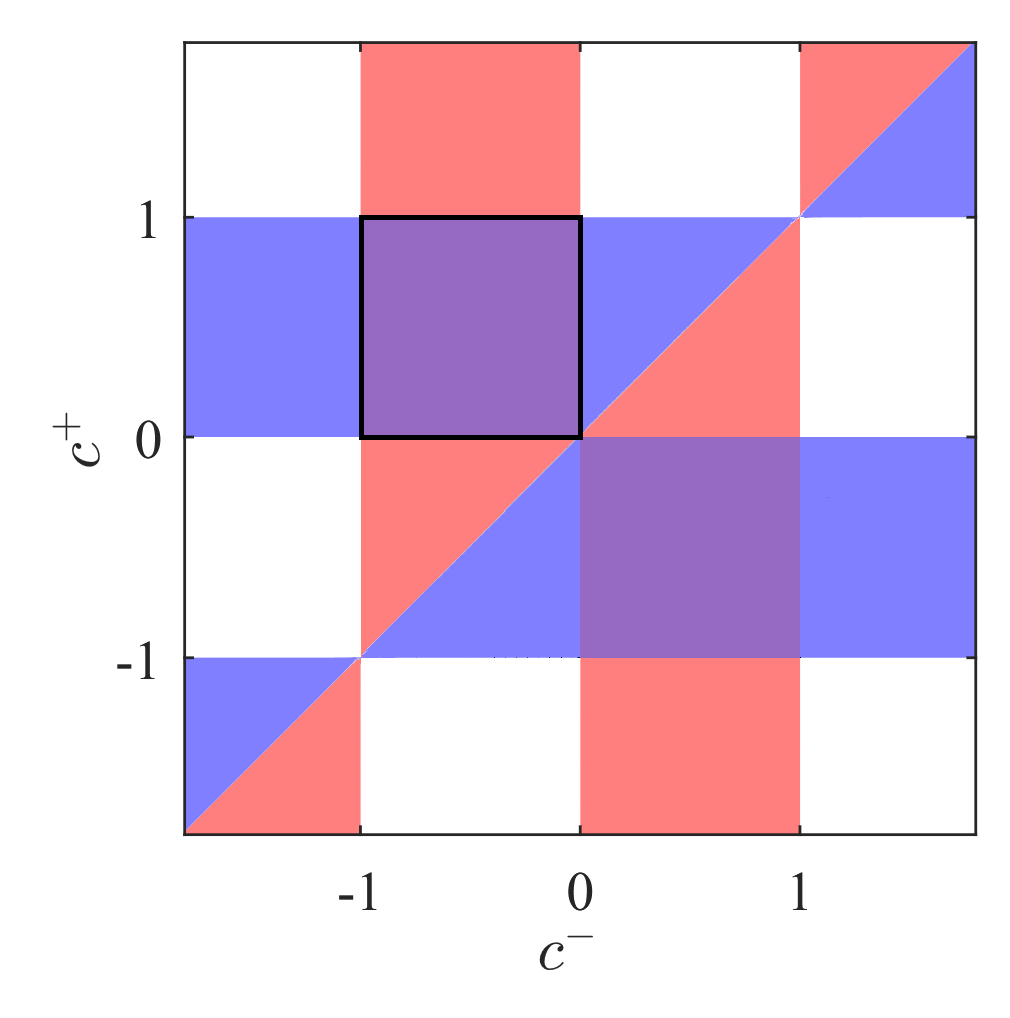}
	\caption{Positivity domain of attenuation rates $\mathcal{R}^{\pm}$ \eqref{eq:attenuation_rates_1D} as a function of eigen-modes $c^{\pm}$. Red: Positivity domain of $\mathcal{R}^{-}$. Blue: Positivity domain of $\mathcal{R}^{+}$. Purple: Positivity domain of both attenuation rates $\mathcal{R}^{\pm}$ simultaneously. The positive square root convention restricts the stability domain to the top-left quadrant, shown with solid black lines.}
	\label{R_positivity}
\end{figure}

Moving on to the dissipation at the NS level, the balance equations become {(see Appendix \ref{app:CE})},
\begin{align}
    &\partial_t \rho + \partial_x (\rho u) = 0, \label{eq:continuity}\\
    &\partial_t (\rho u) + \partial_x \left(\rho u^2+\rho\pi^*+\rho\pi^{\rm neq}\right)= 0,\label{eq:momentum}
\end{align}
where the non-equilibrium pressure be written as,
\begin{equation}\label{eq:noneq_second_order_moment}
    \rho\pi^{\rm neq} = -\left(\frac{1-\beta}{2\beta}\right) \left(2 \mathcal{A} \rho \varsigma^2 \partial_x u + \mathcal{B} \partial_x \rho\right).
\end{equation}
Here $\mathcal{A}$ ({viscosity factor}) and $\mathcal{B}$ ({compressibility error}) are,
\begin{align}
		&\mathcal{A}=\frac{1}{2\varsigma^2}\left(3\varsigma^2-(c^+)^2-(c^-)^2 - c^+ c^- \right), 
		\label{eq:Ac}\\
			&\mathcal{B}= \frac{3\varsigma^2}{2}\left(c^+ + c^-\right) - \frac{1}{2}\left(({c^+})^3 + ({c^-})^3\right).
		\label{eq:Bc}
	\end{align}
Upon linearization, spectral analysis of the balance equations \eqref{eq:continuity} and \eqref{eq:momentum} demonstrates the following dissipation-dispersion relations for the eigen-modes,
    \begin{align}
    	\omega^{\pm} &= c^\pm k + \mathrm{i} \varsigma^2 \left(\frac{1-\beta}{2\beta}\right)\mathcal{R}^{\pm} k^2 + O\left(k^3\right),
    	\label{eq:omega1D}
    \end{align}
where $k$ is the wave vector, $\mathrm{i}=\sqrt{-1}$, and $\mathcal{R}^{\pm}$ are attenuation rates, written in terms of eigen-modes,
\begin{equation}
	\mathcal{R}^{\pm}=\pm \frac{c^{\pm}\left(3 \varsigma^2 - (c^{\pm})^2\right)}{\varsigma^2(c^+-c^-)}.\label{eq:attenuation_rates_1D}
\end{equation}
Positivity of both attenuation rates simultaneously, $\mathcal{R}^{\pm}\ge 0$, is the necessary condition for the linear stability of the LBGK in the hydrodynamic limit $k\to 0$. Positivity domain of both attenuation rates $\mathcal{R}^{\pm}\ge 0$ is shown in Fig.~\ref{R_positivity}. Fig.~\ref{R_positivity} reveals that the linear stability range in terms of eigen-modes reads, 
\begin{equation}
   0  \leq  c^+ \leq 1,\
     -1\leq  c^- \leq 0. 
     \label{eq:condition1}
\end{equation}
In the expanded dimensional form, condition \eqref{eq:condition1} is written, 
\begin{align}
    0  \leq  u+\varsigma_\rho \leq \frac{\delta x}{\delta t},\
     \frac{\delta x}{\delta t}\leq  u-\varsigma_\rho \leq 0,\label{eq:condition1_SWE}
\end{align}
{Thus, the \emph{necessary} stability condition is established for LBGK with a generic non-ideal equation of state. It should be noted that, in terms of the eigen-modes $c^\pm$, condition \eqref{eq:condition1} is the same as the one obtained for the ideal equation of state in \cite{PhysRevE.110.015306}. Similar to \cite{PhysRevE.110.015306}, we refer to unconditional linear stability (in the hydrodynamic limit) if condition \eqref{eq:condition1} holds in the entire range of flow velocity $u\in[-1,1]$. Unconditional stability provides for independence on the kinematic viscosity and thus is important for high Reynolds number simulations. However, at variance with \cite{PhysRevE.110.015306}, in the present non-ideal case the density does not cancel out from the eigen-modes \eqref{eq:hydro_modes_D1Q3} but rather contributes via the density-dependent speed of sound \eqref{eq:sound}.}

\begin{figure}[h]
	\centering		\includegraphics[width=0.7\linewidth]{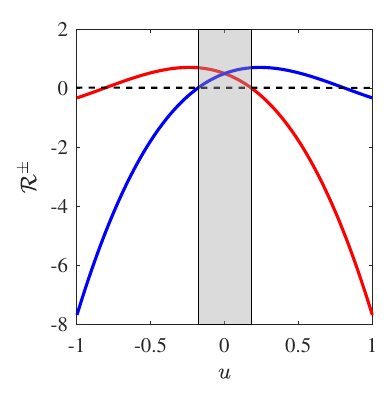}
	\caption{Attenuation rates of normal modes \eqref{eq:attenuation_rates_1D}, $\mathcal{R}^+$ in red and $\mathcal{R}^-$ in blue, for shallow water equation of state \eqref{eq:swe} for ${\rho}=1$ and $g=2/3$. Shaded area represents the stability domain $u\in [-0.1835, 0.1835]$.}
	\label{Rpm_rate_swe}
\end{figure}

\begin{figure}[h]
	\centering		\includegraphics[width=0.7\linewidth]{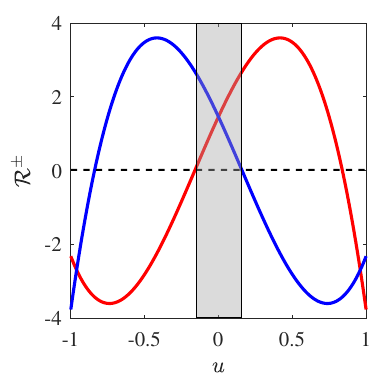}
	\caption{Attenuation rates of normal modes \eqref{eq:attenuation_rates_1D}, $\mathcal{R}^+$ in red and $\mathcal{R}^-$ in blue, for van der Waals equation of state \eqref{eq:vdw} for $T_r=0.8$, $a=1/49$, $b=2/21$, $\rho_r=0.24$. Shaded area represents stability domain $u\in [-0.1605, 0.1605]$.}
	\label{Rpm_rate_vdw}
\end{figure}

{Let us consider two typical examples. Speed of sound \eqref{eq:sound} for the shallow water pressure \eqref{eq:swe} and the van der Waals pressure \eqref{eq:vdw} are, respectively,
\begin{align}\label{eq:swe_drho}
&\varsigma_{\rho} = \sqrt{g \rho},\\
\label{eq:vdw_drho}
&\varsigma_{\rho} = \sqrt{\frac{R T}{{(1-b\rho)}^2} - 2a\rho}.
\end{align}
Dissipation rates  for the shallow water \eqref{eq:swe} and for van der Waals \eqref{eq:vdw} equations of state are shown in Fig.\ \ref{Rpm_rate_swe} and Fig.\ \ref{Rpm_rate_vdw}, respectively.}
{
Stability domain, shown in gray in Figs.~\ref{Rpm_rate_swe} and \ref{Rpm_rate_vdw}, corresponds to the range of the flow velocity where both $\mathcal{R}^+$ and $\mathcal{R}^-$ are positive. In none of these typical cases the necessary stability condition  \eqref{eq:condition1} holds unconditionally, and stable operation range with respect to flow velocity is restricted.}

\subsection{Sufficient stability condition: Arbitrary lattice Knudsen number\label{sec:stability_D1Q3}}

{While the above stability condition \eqref{eq:condition1} was established in the low Knudsen number limit $k\to 0$, we next extend the linear stability analysis to the full discrete lattice Boltzmann system \eqref{eq:LBGK}. To that end, we follow the algorithm of Schur-stability analysis already presented for the LBM in \cite{PhysRevE.110.015306}. The analysis is straightforward but tedious, and we only present the final results.}

\begin{figure}[h]
	\centering		\includegraphics[width=0.9\linewidth]{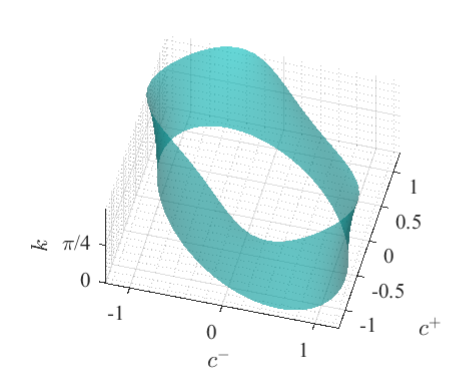}
	\caption{Condition \eqref{eq:cond_3} in the $(c^-,c^+,k)$ parameter space. {By periodicity of the function $\cos^2(k/2)$, only the interval $0\le k\le \pi/2$ is shown.} The interior of the iso-surface represents the inequality \eqref{eq:cond_3}. }
	\label{stability_3d}
\end{figure}

\begin{figure}[h]
	\centering		\includegraphics[width=0.7\linewidth]{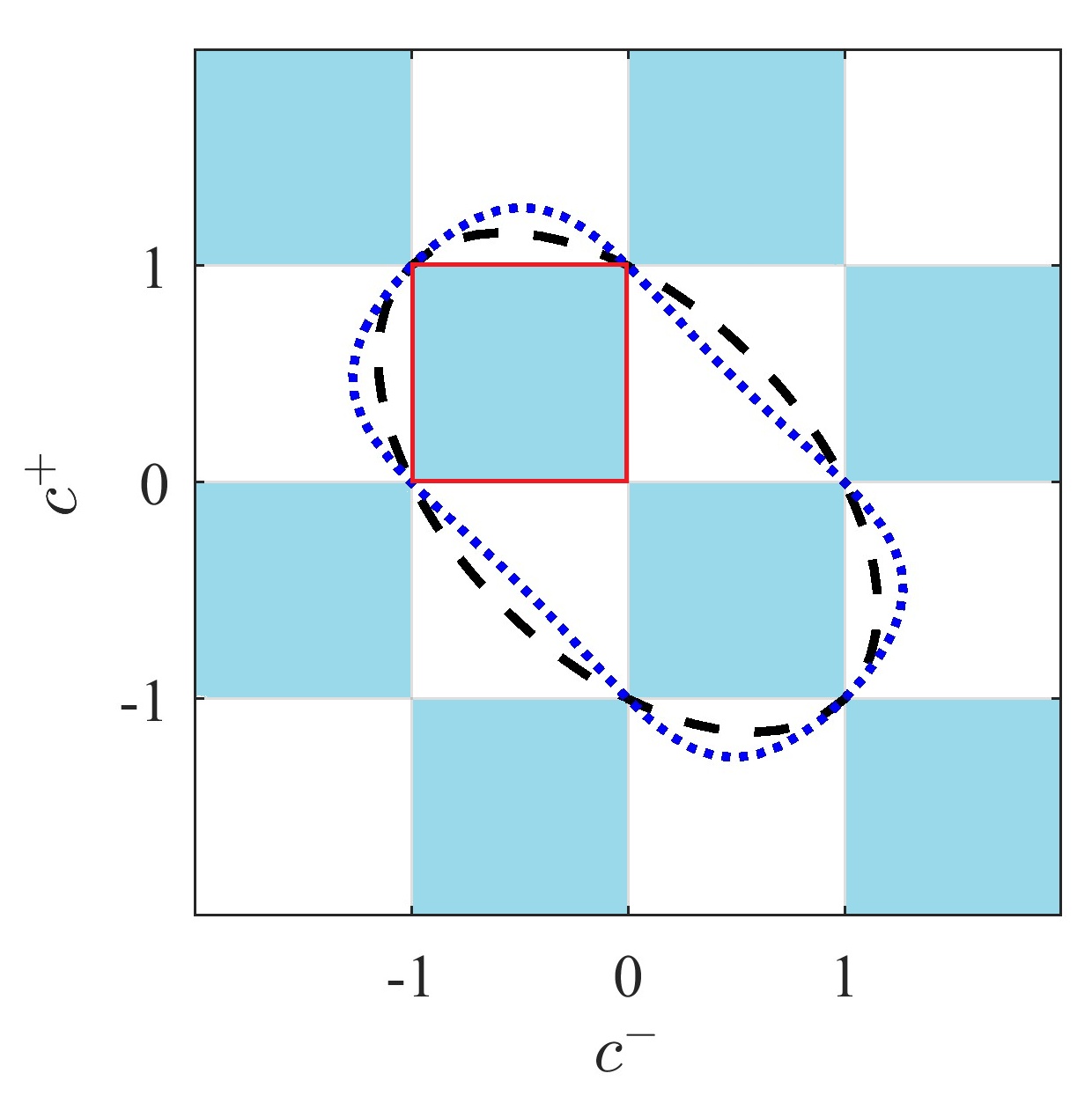}
	\caption{{Linear stability conditions of the D$1$Q$3$ LBGK system. Blue-colored areas: Condition \eqref{eq:cond_2}. Intersection of the interiors of the dashed ($k=0$) and of the dotted ($k=\pi/2$) contours: Condition \eqref{eq:cond_3} for all wave-numbers $k\in[0,2\pi]$ by symmetry. Intersection of the latter interiors with the checkerboard colored area marked with red contour represents the domain of linear stability of the D$1$Q$3$ LBGK system.}}
	\label{D1Q3_linear_stability}
\end{figure}

{The most general linear stability condition for the D$1$Q$3$ LBGK system \eqref{eq:LBGK} for the entire range of wave-number, requires that the following three inequalities hold simultaneously:}
\begin{widetext}
\begin{align}
&\beta\left(\beta-1\right){\leq 0},\label{eq:cond_1}\\
&c^+c^-\left(({c^+})^2-1\right) \left(({c^-})^2-1\right)\leq0.\label{eq:cond_2}\\
&c^+c^- - \cos^2\left(\frac{k}{2}\right) + \left(({c^-})^2+({c^+})^2\right)\cos^2\left(\frac{k}{2}\right)+ ({c^-})^2 ({c^+})^2\left(1-\cos^2\left(\frac{k}{2}\right)\right)\leq0. \label{eq:cond_3}
\end{align}
\end{widetext}
{Condition \eqref{eq:cond_1} is decoupled from the two other conditions and is always satisfied since $\beta\in[0,1]$. The two remaining conditions \eqref{eq:cond_2} and \eqref{eq:cond_3} are formulated in terms of the hydrodynamic eigen-modes $c^\pm$ \eqref{eq:modes_D1Q3_gen} and the wave-number $k$. Dependence on the wave-number is manifest only in the third condition \eqref{eq:cond_3} which represents a family of symmetric quartic curves, parameterized by $k$, with an ellipse in the limit $k\to 0$. Condition \eqref{eq:cond_3} is shown in the $(c^-,c^+,k)$ parameter space in Fig.~\ref{stability_3d}.}
{Inequality \eqref{eq:cond_2} holds in a checkerboard domain in the $(c^-,c^+)$ plane shown in Fig.\ \ref{D1Q3_linear_stability}, and is independent of the wave-number $k$.}  
{Moreover, the domain of eigen-modes where both linear stability conditions \eqref{eq:cond_2} and  \eqref{eq:cond_3} are satisfied simultaneously for all wave-numbers $k$ is also shown in  Fig.~\ref{D1Q3_linear_stability}.}
{Comparison of the stability domains in Figs.\ \ref{D1Q3_linear_stability} and \ref{R_positivity} demonstrates that, the hydrodynamic limit stability condition \eqref{eq:condition1} is at the same time also the sufficient linear stability condition for all Knudsen numbers.
This confirms that, very much like the analysis presented in \cite{PhysRevE.110.015306} for the ideal gas pressure, the condition derived in the hydrodynamic limit holds for all values of the wave-number. This completes the detailed analysis of the D$1$Q$3$ lattice resulting in the stability condition \eqref{eq:condition1} in terms of the normal eigen-modes propagation speeds. We next proceed with a linear stability analysis of the two-dimensional LBGK system.}
\section{Stability of the two-dimensional system\label{sec:D2Q9}}
\subsection{Spectral analysis of the target equations}
The target continuity and momentum Navier--Stokes equations are,
	\begin{align}
  &  \partial_t \rho + \bm{\nabla}\cdot\rho\bm{u} = 0,\label{eq:target_rho}\\
   & \partial_t \rho \bm{u} + \bm{\nabla}\cdot\rho\bm{u}\otimes\bm{u} + \bm{\nabla}\cdot{\rho}\pi^*\bm{I} + \bm{\nabla}\cdot\bm{T} = 0,\label{eq:target_u}
	  \end{align}
where the viscous stress tensor is defined as,
\begin{equation}\label{eq:NS_shear}
    \bm{T} = -\rho \nu\left(\bm{\nabla}\bm{u}+\bm{\nabla}\bm{u^\dagger} - \bm{\nabla}\cdot\bm{u}\bm{I}\right) - {\rho}\eta \bm{\nabla}\cdot\bm{u}\bm{I},
\end{equation}
and where $\nu$ and $\eta$ are the shear and the bulk viscosity, respectively.
{Spectral analysis of the linearized target equations reveals the following dispersion-dissipation relations in two dimensions, see Appendix \ref{app:NS_non_ideal}:  
\begin{align}
   & \omega_{\rm shear} = \bm{k}\cdot\bm{u} + \mathrm{i}\nu {k}^2,\label{eq:shear_mode_NS}\\
   & \omega_{{\rm ac}\pm} = \bm{k}\cdot\bm{u} \pm  {k}\varsigma_\rho + \mathrm{i}\frac{(\nu+\eta)}{2} {k}^2.\label{eq:normal_mode_NS}
	  \end{align}
Here $\bm{k}$ the wave-vector and $k=|\bm{k}|$ the wave-number.
The dispersion and dissipation of the D$2$Q$9$ LBGK model shall be probed against these target equations.}
\subsection{Hydrodynamic limit of the LBGK model} 
{The standard Chapman--Enskog analysis of the LBGK system \eqref{eq:LBGK} recovers the following balance equations, see Appendix \ref{app:CE}:
	\begin{align}
   & \partial_t \rho + \bm{\nabla}\cdot\rho\bm{u} = 0,\label{eq:cons_0_CE}\\
   & \partial_t \rho \bm{u} + \bm{\nabla}\cdot\rho\bm{u}\otimes\bm{u} + \bm{\nabla}\cdot{\rho}\pi^*\bm{I} + \bm{\nabla}\cdot{\rho}\bm{\pi}^{\rm neq} = 0,\label{eq:cons_1_CE}
	  \end{align}
where the non-equilibrium pressure tensor is written as,
\begin{equation}\label{eq:neq_2d}
   {\rho} \bm{\pi}^{\rm neq} = \bm{T} + \bm{T}'.
\end{equation}
{Here $\bm{T}$ is the target viscous stress tensor \eqref{eq:NS_shear} while the deviation $\bm{T}'$ has the form,}
\begin{equation}\label{eq:Tprime}
    \bm{T}' = - \left(\frac{1-\beta}{2\beta}\right)\left[ \rho \varsigma^2\bm{A}'\circ\left(\bm{I} \circ \bm{\nabla}\bm{u}\right) +  \bm{B}\circ\left(\bm{I} \circ {\rm diag}(\bm{\nabla}\rho)\right)\right],
\end{equation}
where $\circ$ denotes the Hadamard product while  $\bm{A}'$ and $\bm{B}$ are diagonal matrices,
\begin{align}
\label{eq:error_A_2d}
&\bm{A}' = {\rm diag}(\mathcal{A}'_{\alpha\alpha}), &   \mathcal{A}'_{\alpha\alpha} = -\frac{3 u_\alpha^2 + 3(\pi^* - \varsigma^2)}{\varsigma^2},\\
\label{eq:error_B_2d}
&\bm{B} = {\rm diag}(\mathcal{B}_{\alpha\alpha}),   &    \mathcal{B}_{\alpha\alpha} = - u_\alpha^3 - 3u_\alpha(\varsigma_\rho^2 - \varsigma^2).
\end{align}
For the D$2$Q$9$ lattice, the deviation \eqref{eq:Tprime} has the following explicit form:
\begin{equation}
    \bm{T}' = -\left(\frac{1-\beta}{2\beta}\right)\begin{bmatrix} \rho \varsigma^2 \mathcal{A}'_{xx} \partial_x u_x + \mathcal{B}_{xx} \partial_x \rho & 0\\ 0 & \rho \varsigma^2 \mathcal{A}'_{yy} \partial_y u_y + \mathcal{B}_{yy} \partial_y \rho\end{bmatrix},
\end{equation}
Finally, the kinematic viscosity $\nu$ is tied to $\beta$ through \eqref{eq:viscosity} while the bulk viscosity $\eta$ reads,
\begin{equation}\label{eq:bulk_visc}
    \eta = \nu \left(\frac{2+D}{D} - \frac{\partial \ln(\rho\pi^*)}{\partial\ln \rho}\right).
\end{equation}
Note that, the deviation  \eqref{eq:Tprime} from the target Navier--Stokes viscous stress tensor \eqref{eq:NS_shear} only affects the diagonal components of the LBGK nonequilibrium pressure tensor \eqref{eq:neq_2d}.}

{In order to quantify the deviations, and to establish relations to stability conditions in the above one-dimensional case, we perform spectral analysis of Eqs.~\eqref{eq:cons_0_CE} and \eqref{eq:cons_1_CE} considering a perturbation of the following form,
\begin{equation}\label{eq:perturbation2D}
    \{\rho',\bm{u}'\} = \{\hat{\rho},\hat{\bm{u}}\}\exp{[\mathrm{i}(k_\alpha x_\alpha}-\omega t)],
\end{equation}
Under the perturbation \eqref{eq:perturbation2D}, the dispersion-dissipation relations for the eigen-modes of the balance equations \eqref{eq:cons_0_CE} and \eqref{eq:cons_1_CE} are found to be,
\begin{align}
   & \omega_{\rm shear} = k_\alpha u_\alpha + \mathrm{i}\nu k_\alpha^2,\label{eq:share2}\\
   & \omega_{{\rm ac}\pm} = k_\alpha (u_\alpha \pm  \varsigma_\rho) + \mathrm{i} \varsigma^2 \left(\frac{1-\beta}{2\beta}\right)\mathcal{R}^{\pm} k_\alpha^2 + \mathcal{O}(k_\alpha^3).
    \label{eq:normal2}
	  \end{align}
Thus, the normal modes \eqref{eq:normal2} correspond exactly to those in Eq.~\eqref{eq:omega1D}, confirming that \eqref{eq:condition1} remains the necessary stability condition of the D$2$Q$9$ LBGK.}
{At the same time, dissipativity of the shear mode \eqref{eq:share2} requires positivity of the pressure via the shear viscosity relation \eqref{eq:viscosity},
\begin{equation}
    \pi^*\ge0.
\end{equation}
Note that positivity of the pressure is different from the existence of the speed of sound \eqref{eq:sound_existence}, and is not present in the one-dimensional setting due to the absence of the shear mode therein.}
\subsection{Spectral analysis of the D$2$Q$9$ LBGK}
{We shall now analyze spectral dispersion (real part of the dispersion-dissipation relations, $\Re(\omega_i)$, $i=\{{\rm shear},{\rm ac}\pm\}$) and spectral dissipation (imaginary part of the dispersion-dissipation relations, $\Im (\omega_i)$, $i=\{{\rm shear},{\rm ac}\pm\}$). 
To that end, we generated three data sets: First set is the target Navier--Stokes spectra \eqref{eq:shear_mode_NS} and \eqref{eq:normal_mode_NS}. The second set is the spectra of the hydrodynamic limit of the LBGK \eqref{eq:share2} and \eqref{eq:normal2}. The third data set was obtained by spectral analysis of the LBGK equation using standard numerical procedures~\cite{Ali_StabilityEntropic,hosseini2020development}.}
{The shallow water and the van der Waals non-ideal pressures, Eqs.~\eqref{eq:swe} and \eqref{eq:vdw}, were used for the comparison.}

\begin{figure}[h]
	\centering
		\includegraphics[width=0.7\linewidth]{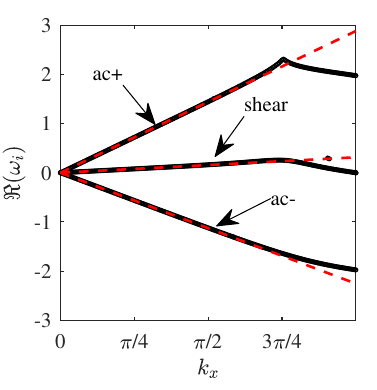}
	\caption{{Spectral dispersion $\Re (\omega_{\rm shear})$ and $\Re(\omega_{{\rm ac}\pm})$ of hydrodynamic modes of the D$2$Q$9$ LBGK for the shallow water pressure \eqref{eq:swe}, with ${u}_x=0.1$, $g=2/3$, $\rho=1$ and $\beta=0.625$. Symbol: Numerical spectral analysis of the LBGK equation \eqref{eq:LBGK};  Red dashed lines: Spectral dispersion of the target Navier--Stokes equations \eqref{eq:shear_mode_NS} and \eqref{eq:normal_mode_NS}.}}
	\label{dispersion_LBM}
\end{figure}

\begin{figure}[h]
	\centering
		\includegraphics[width=0.7\linewidth]{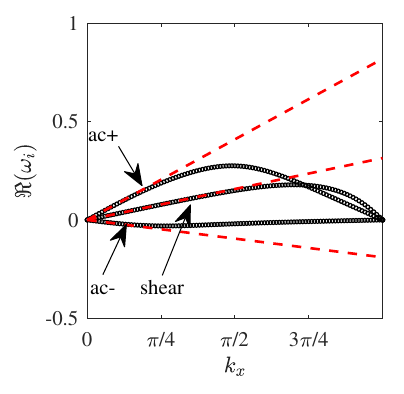}
	\caption{{Spectral dispersion $\Re (\omega_{\rm shear})$ and $\Re(\omega_{{\rm ac}\pm})$ of hydrodynamic modes of the D$2$Q$9$ LBGK for the van der Waals pressure \eqref{eq:vdw}, with ${u}_x=0.1$, $T_r=0.8$, $a=1/49$, $b=2/21$, $\rho_r=0.24$ and $\beta=0.625$.  Symbol: Numerical spectral analysis of the LBGK equation \eqref{eq:LBGK};  Red dashed lines: Spectral dispersion of the target Navier--Stokes equations \eqref{eq:shear_mode_NS} and \eqref{eq:normal_mode_NS}.}}
	\label{dispersion_LBM_vdw}
\end{figure}

{The spectral dispersion plots for these two equations of state are shown in Figs.~\ref{dispersion_LBM} and \ref{dispersion_LBM_vdw}.}
{In the hydrodynamic limit  $k\rightarrow 0$, in agreement with the results of the multi-scale analysis, both plots show convergence of hydrodynamic modes to the spectral dispersions of Eqs.~\eqref{eq:shear_mode_NS} and \eqref{eq:normal_mode_NS}.}

{The spectral dissipation rates of the full LBGK system, obtained as $\Im(\omega_i)$, are shown in Figs.~\ref{dissipation_LBM} and \ref{dissipation_LBM_vdW} for Eqs.~\eqref{eq:swe} and \eqref{eq:vdw}.
\begin{figure}[h]
	\centering
		\includegraphics[width=0.7\linewidth]{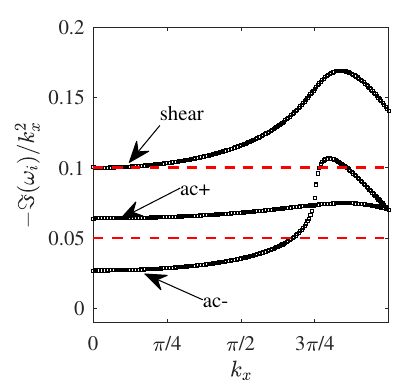}
	\caption{{Spectral dissipation of hydrodynamic modes of the D$2$Q$9$ LBGK for the shallow water equation of state \eqref{eq:swe}, ${u}_x=0.1$, ${\rho}=1$, $g=2/3$ and $\beta=0.625$.
    Symbol: Numerical spectral analysis of the LBGK equation \eqref{eq:LBGK}; Red dashed lines: Spectral dissipation, $-\Im(\omega_{\rm shear})/k^2$ and $-\Im(\omega_{\rm ac^{\pm}})/k^2$, of the target Navier--Stokes equations \eqref{eq:shear_mode_NS} and \eqref{eq:normal_mode_NS} with shallow water equation of state.}}
	\label{dissipation_LBM}
\end{figure}
\begin{figure}[h]
	\centering
		\includegraphics[width=0.7\linewidth]{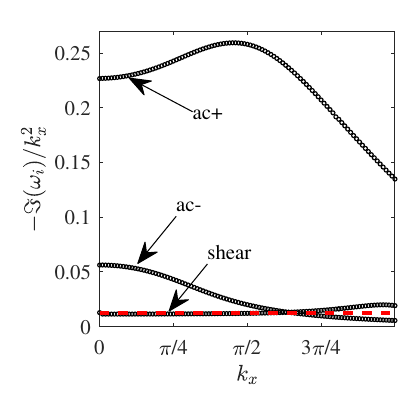}
	\caption{{Spectral dissipation of hydrodynamic modes of the D$2$Q$9$ LBGK for the van der Waals equation of state \eqref{eq:vdw}, ${u}_x=0.1$, $T_r=0.8$, $a=1/49$, $b=2/21$, $\rho_r=0.24$ and $\beta=0.625$.
    Symbol: Numerical spectral analysis of the LBGK equation \eqref{eq:LBGK}; Red dashed lines: Spectral dissipation, $-\Im(\omega_{\rm shear})/k^2$ and $-\Im(\omega_{\rm ac^{\pm}})/k^2$, of the target Navier--Stokes equations \eqref{eq:shear_mode_NS} and \eqref{eq:normal_mode_NS} with shallow water equation of state.}}
	\label{dissipation_LBM_vdW}
\end{figure}
Both plots show that in the hydrodynamic limit, while the shear mode dissipation rate converges to that of the target Navier--Stokes equations, the normal modes show deviations. The deviations of the dissipation rates are shown in the entire range of velocity in Figs.~\ref{hydro_dissipation_swe} and \ref{hydro_dissipation_vdw}.
\begin{figure}[h]
	\centering
		\includegraphics[width=\linewidth]{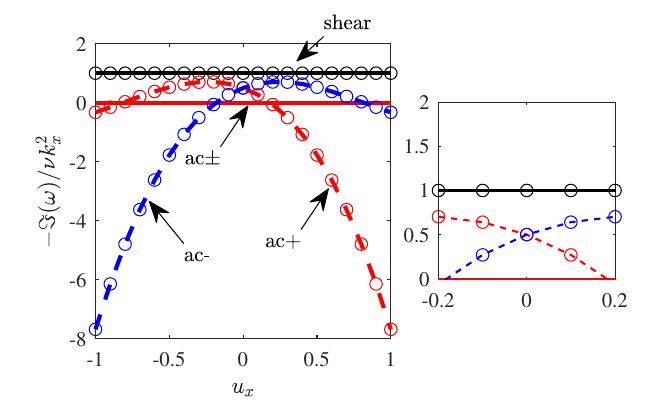}
	\caption{{Spectral dissipation of hydrodynamic modes in the limit of $k\rightarrow 0$ for the shallow water equation \eqref{eq:swe} with $g=2/3$, $\rho=1$ and $\beta=0.625$. Solid lines: Eqs.~\eqref{eq:shear_mode_NS} and \eqref{eq:normal_mode_NS}, dashed lines: Eqs.~\eqref{eq:share2} and \eqref{eq:normal2}, circular markers: data from the full LBGK spectral analysis. Inset: Zoomed in region around $u=0$.}}
	\label{hydro_dissipation_swe}
\end{figure}
\begin{figure}[h]
	\centering
		\includegraphics[width=\linewidth]{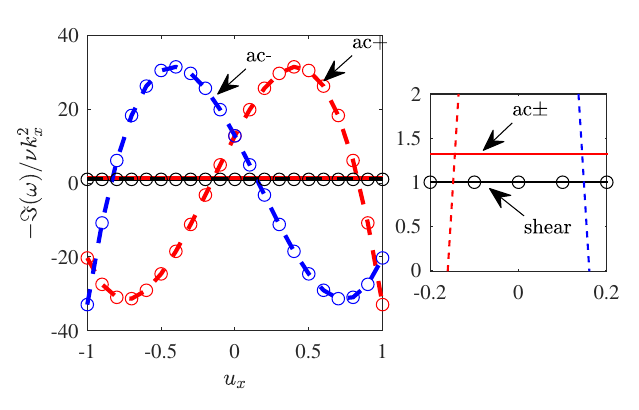}
	\caption{{Spectral dissipation of hydrodynamic modes in the limit of $k\rightarrow 0$ for the van der Waals equation \eqref{eq:vdw} with $T_r=0.8$, $a=1/49$, $b=2/21$, $\rho_r=0.24$ and $\beta=0.625$. Solid lines: Eqs.~\eqref{eq:shear_mode_NS} and \eqref{eq:normal_mode_NS}, dashed lines: Eqs.~\eqref{eq:share2} and \eqref{eq:normal2}, circular markers: data from the full LBGK spectral analysis. Inset: Zoomed in region around $u=0$.}}
	\label{hydro_dissipation_vdw}
\end{figure}
In the limit $k\rightarrow 0$, the correct shear mode dissipation rate is recovered. However, the normal modes' dissipation rates shows deviations. These deviations originate from the lattice constraint, $c^3_{i\alpha}=c_{i\alpha}$, which leads to the biased third-order moments of the D$D$Q$3^D$ lattices.
While not the focus of the present work, we mention that the deviations in the dissipation rate  of normal modes at $k\to 0$ can be removed via corrections to the the non-equilibrium second-order moment, see \cite{prasianakis2009lattice,li2012coupling,hosseini2020compressibility} for ideal fluids and \cite{hosseini2022towards} for non-ideal equations of state.
}
\subsection{Linear stability domain of D$2$Q$9$ LBGK}
{Unlike the fully analytical study of linear stability of the one-dimensional D$1$Q$3$ LBGK in sec.\ \ref{sec:stability_D1Q3} resulting in the analytical stability conditions \eqref{eq:cond_2} and \eqref{eq:cond_3}, Schur-stability analysis becomes formidable in a higher dimension.}
Thus, we restore to a numerical study of the linear stability domain of the D$2$Q$9$ LBGK. We consider the parameter space $(\varsigma_\rho, \bm{u}, \beta)$. This choice is motivated by the fact that the analytical stability condition of the D$1$Q$3$ lattice \eqref{eq:attenuation_rates_1D} can be equivalently written as a function of $\varsigma_\rho$ and $u$ using Eq.~\eqref{eq:hydro_modes_D1Q3},
\begin{align}\label{eq:connection_1} 
 &  c^+ + c^- = 2u,\\
&  \left(c^+-u\right) \left(c^- -u\right) = -\varsigma^2_\rho.\label{eq:connection_2}
\end{align}
Stability was studied following the protocol described in \cite{hosseini2020development}. The maximal flow velocity magnitude was computed for which the growth rate $\Im{(\omega_i)}$ {of all modes of the linearized LBGK equation $\omega_i$, $i={1,\dots,9}$,} remains negative  for $(k_x,k_y)\in [-\pi,\pi]\times[0,\pi]$ and angles of the flow velocity $\bm{u}$ relative to the $x$-axis in $[-\pi/2,\pi/2]$. The wave-number interval was spaced with $\Delta k = 0.02$, and the velocity vector angles with the $x$-axis was spanned in steps of $\pi/20$. 

\begin{figure}[h]
	\centering
		\includegraphics[width=0.8\linewidth]{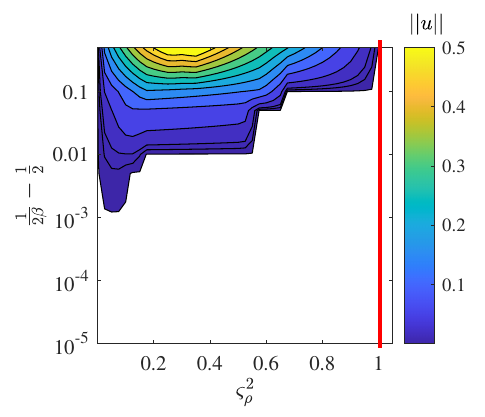}
	\caption{Linear stability domain for D2Q9 LBGK with the shallow water equation of state \eqref{eq:swe}, $g=2/3$.
    {Red vertical line at $\varsigma^2_\rho=1$: the necessary stability condition of normal modes \eqref{eq:condition1_SWE} at $u=0$. No stable realizations are possible to the right of the red line. Complex pattern of the stability domain on the left of the red line demonstrates the lack of unconditional stability.}}
	\label{stability_2d_swe}
\end{figure}

Linear stability domain  for the shallow water pressure \eqref{eq:swe} is shown in Fig.~\ref{stability_2d_swe}.
{For the van der Waals equation, we first computed the saturated liquid densities for reduced temperatures $T/T_c\in[0.4,0.91]$. These liquid density values, along with the corresponding temperatures, were used to establish the pressure \eqref{eq:vdw} and the speed of sound \eqref{eq:vdw_drho} used in the linear stability analysis.}
\begin{figure}[h]
	\centering
	\includegraphics[width=0.8\linewidth]{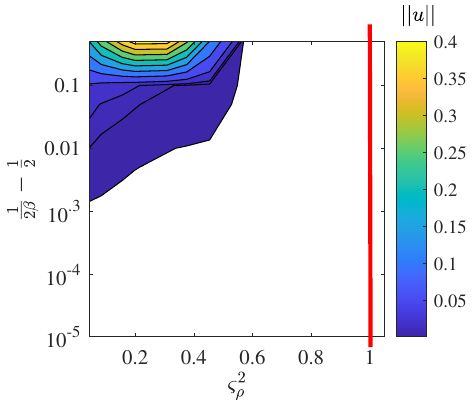}
	\caption{Linear stability domain for D$2$Q$9$ LBGK with the van der Waals equation of state \eqref{eq:vdw}, with $a=1/49$, $b=2/21$, $R=1$ and $T/T_c\in[0.4,0.91]$. Densities on the liquid branch of the saturation line are used.
     {Red vertical line: Same as in Fig.\ \ref{stability_2d_swe}.}}
	\label{stability_2d_vdw}
\end{figure}
Results of the stability analysis for the van der Waals pressure are reported in Fig.~\ref{stability_2d_vdw}.

\begin{figure}[h]
	\centering
	\includegraphics[width=0.8\linewidth]{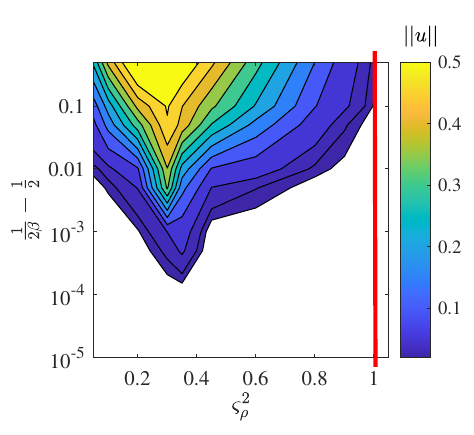}
	\caption{Linear stability domain for D$2$Q$9$ LBGK with the ideal gas equation of state \eqref{eq:ideal_eos}, with $R=1$ and $T\in[0\,1]$.
    {Red vertical line: Same as in Fig.\ \ref{stability_2d_swe}.}
    }
	\label{stability_LBM_ideal}
\end{figure}

First, we note that the linear stability condition \eqref{eq:condition1} derived for the D$1$Q$3$ system still holds in higher dimension but only as a necessary (rather than sufficient) condition. This is best seen  by that all stability plots show unconditional \emph{instability} for $\varsigma_\rho\geq1$.
The linear stability domains in Figs.\  \ref{stability_2d_swe} and \ref{stability_2d_vdw} show qualitative similarities with the ideal equation of state,
\begin{equation}\label{eq:ideal_eos}
    \pi^* = R T,
\end{equation}
shown in Fig.~\ref{stability_LBM_ideal}: Maximal flow speeds are achieved around $\varsigma_\rho=\varsigma$ while the interval $0\leq\varsigma_\rho\leq\varsigma$ allows for stable simulations at higher $\beta$ as compared to $\varsigma_\rho \geq\varsigma$, where LBGK becomes unconditionally unstable. For example, in Fig.~\ref{stability_2d_vdw}, instability sets in for $\varsigma_\rho > 0.74$.
At variance with the ideal equation of state, even at $0\leq\varsigma_\rho\leq\varsigma$, the solvers are unconditionally unstable for $\beta\geq 0.98$.
{We remind that, for the speed of sound \eqref{eq:sound} equal to the lattice speed of sound \eqref{eq:cs_lattice}, $\varsigma_\rho=\varsigma=1/\sqrt{3}$, the \emph{absolute} maximum of the flow velocity $u^{\rm max}$ above which no linear stability of the lattice Boltzmann models is possible on the D$D$Q$3^D$ lattice with the pressure of the form \eqref{eq:Peq_gen} is \cite{Ali_StabilityEntropic,PhysRevE.110.015306},
\begin{equation}
    \label{eq:umax}
    u^{\rm max}=\left(1-\frac{1}{\sqrt{3}}\right)\approx 0.4226,
\end{equation}
This is consistent with all stability plots as no stability is observed for flow velocities higher than $u^{\rm max}$ for $\varsigma_\rho\geq \varsigma$.}
A discussion is presented in the next section. 

\section{{Discussion and outlook}\label{sec:conslusions}}
{We presented a systematic linear stability analysis of a class of lattice Boltzmann models with non-ideal equations of state. Following the methodology set out in our previous publications \cite{Ali_StabilityEntropic,PhysRevE.110.015306}, we provided a comprehensive linear stability analysis of the one-dimensional D$1$Q$3$ LBGK system.
The linear stability condition was found as a bound on normal eigen-modes \eqref{eq:condition1}, which confirms that the previous result for ideal gas equation of state \cite{PhysRevE.110.015306} also holds for the general non-ideal pressure. Remarkably, the full linear stability of the D$1$Q$3$ LBGK model is solely defined by the hydrodynamic eigen-modes \eqref{eq:modes_D1Q3_gen} and \eqref{eq:sound} through their boundedness \eqref{eq:condition1}. Geometrical feature of the linear stability condition was highlighted in Fig.~\ref{D1Q3_linear_stability}.
These general findings were exemplified by two typical nonlinear equations of state, the shallow water and the van der Waals. 
It was demonstrated that conventional equations of state and the standard equilibrium pressure \eqref{eq:Peq_gen} do not satisfy the stability requirement \eqref{eq:condition1} unconditionally for the entire range of flow velocity.}

{The analysis was extended to higher dimensions, where stability conditions of the one-dimensional system were shown to remain necessary. Spectral dissipation and dispersion properties were studied for the two-dimensional D$2$Q$9$ LBGK model. 
Spectral analysis revealed that LBGK  captures correctly the shear mode while admitting deviations in the normal modes dissipation rates.
Finally, the linear stability domain was studied numerically for the D$2$Q$9$ LBGK model. The numerical analysis confirmed once again that boundedness of the normal modes \eqref{eq:condition1} remains the necessary (but not sufficient) stability condition. 
Overall, it was observed that the linear stability range for typical non-ideal equations of state was limited even more than for the ideal gas equation of state.}

\begin{figure}[h!]
	\centering
	\includegraphics[width=6cm,keepaspectratio]{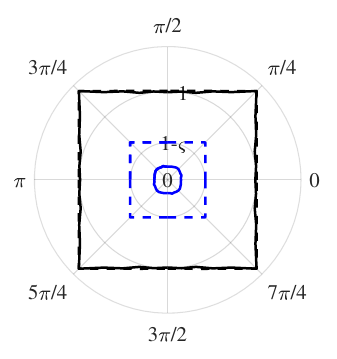}
	\caption{Directional stability domain of the isothermal D$2$Q$9$ LBGK model with (black) the  entropic \eqref{eq:Pstar_entropic} and (blue) standard pressure, $\pi^*=\varsigma^2$, for two different values of kinematic viscosity: (low viscosity, dashed) $\nu=10^{-5}$ and (high viscosity, solid) $\nu=0.1$. {Unconditional linear stability of the entropic LBGK: Stability is manifest on the entire domain of the flow velocity, $\bm{u}\in[-1,1]\times[-1,1]$ independently of the viscosity.}}
	\label{Fig:stability_all}
\end{figure}

{In the outlook, we comment that present study suggests a route to extend the linear stability domain via construction of special equations of state allowing for a non-trivial dependence on the flow velocity, $\pi^*=\pi^*(\rho,u)$.
To that end, we remind the example of isothermal pressure $\pi^*=\varsigma^2$. In this case, the entropic equilibrium with,
\begin{align}
	\pi^{*}& = \varsigma^2 \left( 2\sqrt{{1+(u/\varsigma)}^2 } -{1-(u/\varsigma)}^2\right),\label{eq:Pstar_entropic}
\end{align}
leads to \emph{unconditional} linear stability, as was demonstrated in in \cite{Ali_StabilityEntropic,PhysRevE.110.015306}.
Linear stability domain of the entropic LBGK is juxtaposed to the LBGK with standard equilibrium in Fig.\ \ref{Fig:stability_all}: Unlike the latter, the entropic version is linearly stable on the entire domain of the flow velocity, $\bm{u}\in[-1,1]\times[-1,1]$, independently of the viscosity. 
Following this example, we recognize two key considerations towards stable LBGK models: consistency and unconditional stability. 
The consistency condition can be written down as,
\begin{equation}\label{eq:consistency_general}
    \pi^*(\rho,u) = \pi^*(\rho)\left[1+\mathcal{O}{\left(\frac{u}{\varsigma_\rho}\right)}^4\right],\ \frac{u}{\varsigma_\rho}\to 0.
\end{equation}
In other words, the assumed dependence of the equation of state on the flow velocity should be weaker than the quadratic term $u_\alpha^2$ in the equilibrium pressure \eqref{eq:Peq_gen}.
Eigen-modes of the D$1$Q$3$ LBGK become, instead of \eqref{eq:hydro_modes_D1Q3},}
\begin{align}\label{eq:modes_D1Q3_gen}
	c^{\pm}= u + \frac{\partial_u \pi^*}{2}
	\pm \sqrt{{\left(\frac{\partial_u \pi^*}{2}\right)}^2 + \partial_\rho (\rho\pi^*)}.
\end{align}
{Note additional terms in this more general expression due to a dependence of the pressure on the flow velocity. Nevertheless, the stability condition \eqref{eq:condition1} still holds for the more general case \eqref{eq:modes_D1Q3_gen}. The above example of the entropic equilibrium \eqref{eq:Pstar_entropic} confirms to both the consistency and the stability conditions.}
{Thus, conditions \eqref{eq:consistency_general} and \eqref{eq:condition1} form a basis for designing unconditionally stable LBGK realizations for nonideal equations of state which will be the topic of our future studies.
}
\section*{Acknowledgments}
{We thank the organization of the 33rd Discrete Simulation of Fluid Dynamics Conference (DSFD) in Zurich, Switzerland held over July 9-12 2024 for creating the platform at which
this work was first presented.}
This work was supported by European Research Council (ERC) Advanced Grant No. 834763-PonD and by the Swiss National Science Foundation (SNSF) Grant 200021-228065. Computational resources at the Swiss National Super Computing Center (CSCS) were provided under Grant No. s1286. 
\section*{Author Declarations}

\noindent\textbf{Conflict of Interest}\\
The authors have no conflicts to disclose.

\noindent\textbf{Ethical approval}\\
The work presented here by the authors did not require ethics approval or consent to participate.

\section*{Appendix}
\appendix
\section{Multi-scale analysis of the generic LBM with non-linear-in-density pressure\label{app:CE}}
Here we present a multi-scale perturbation analysis of the LBGK equation \eqref{eq:LBGK}. We first perform a Taylor series expansion around $(\bm{r},t)$:
\begin{equation}\label{eq:taylor}
	    \delta t \mathcal{D}_t f_i + \frac{{\delta t}^2}{2}{\mathcal{D}_t}^2 f_i + {O}({\delta t}^3)= 2\beta\left(f_i^{\rm eq} - f_i\right),
\end{equation}
where we have retained terms up to order two and where $\mathcal{D}_t = \partial_t + \bm{c}_i\cdot\bm{\nabla}$. Introducing characteristic flow scale $\mathcal{L}$ and characteristic flow velocity $\mathcal{U}$, equation \eqref{eq:taylor} is made non-dimensional,
\begin{equation}
	   \left(\frac{\delta r}{\mathcal{L}}\right) \mathcal{D}_t' f_i + \frac{1}{2}{\left(\frac{\delta r}{\mathcal{L}}\right)}^2{\mathcal{D}_t'}^2 f_i = 2\beta\left(f_i^{\rm eq} - f_i\right),
\end{equation}
where primed variables denote non-dimensional form and
\begin{equation}
	\mathcal{D}_t' = \frac{\mathcal{U}}{c}\left(\partial_t' + \bm{c}_i\cdot\bm{\nabla}'\right),
\end{equation}
where $c=\delta r/\delta t$. Assuming acoustic, i.e. $\frac{\mathcal{U}}{c}\sim 1$ and hydrodynamic, i.e.
\begin{equation}\label{eq:epsilon}
    \frac{\delta r}{\mathcal{L}}\sim\varepsilon,
\end{equation}
scaling and dropping the primes for the sake of readability:
\begin{equation}
	\varepsilon \mathcal{D}_t f_i + \frac{1}{2}\varepsilon^2{\mathcal{D}_t}^2 f_i + {O}(\varepsilon^3)= 2\beta\left(f_i^{\rm eq} - f_i\right).
\end{equation}
Then introducing multi-scale expansions:
\begin{subequations}
    \begin{align}
        f_i &= f_i^{(0)} + \varepsilon f_i^{(1)} + \varepsilon^2 f_i^{(2)} + O(\varepsilon^3),\\
        \partial_t &= \varepsilon \partial_t^{(1)} + \varepsilon^2 \partial_t^{(2)} + O(\varepsilon^3),\\
        \bm{\nabla} &= \varepsilon \bm{\nabla},
    \end{align}
\end{subequations}
the following equations are recovered at scales $\varepsilon$ and $\varepsilon^2$:
\begin{subequations}
	\begin{align}
		\varepsilon &: \mathcal{D}_{t}^{(1)} f_i^{(0)} = -2\beta f_i^{(1)},\\
		\varepsilon^2 &: \partial_t^{(2)}f_i^{(0)} + \mathcal{D}_{t}^{(1)} \left(1-\beta\right)f_i^{(1)} = -2\beta f_i^{(2)},
	\end{align}
	\label{Eq:CE_Eq_orders}
\end{subequations}
with $f_i^{(0)}=f_i^{\rm eq}$. Taking the moments of the Chapman-Enskog-expanded equation at order $\varepsilon$:
\begin{subequations}
    \begin{align}
	\partial_t^{(1)}\rho + \bm{\nabla}\cdot\rho \bm{u} &= 0,\label{eq:approach2_continuity1_app}\\
	\partial_t^{(1)}\rho \bm{u} + \bm{\nabla}\cdot\rho \bm{u}\otimes\bm{u} + \bm{\nabla}\cdot\rho\pi^*\bm{I} &= 0.\label{eq:approach2_NS1}
    \end{align}
\end{subequations}
At order $\varepsilon^2$ the continuity equation is:
\begin{equation}
	\partial_t^{(2)}\rho = 0.\label{eq:approach2_continuity2}
\end{equation} 
Summing up Eqs.~\eqref{eq:approach2_continuity1_app} and \eqref{eq:approach2_continuity2} we recover the continuity equation as:
\begin{equation}
	\partial_t \rho + \bm{\nabla}\cdot\rho\bm{u} + \mathcal{O} (\varepsilon^3) = 0.
\end{equation}
For the momentum balance equation, using the mass and momentum balance equations from order $\varepsilon$ and,
\begin{equation}
    \partial_t^{(1)}\rho\pi^* + \bm{\nabla}\cdot\rho\pi^*\bm{u} = \left(\rho\pi^*-\rho\partial_\rho \rho\pi^*\right)\bm{\nabla}\cdot\bm{u},
\end{equation}
we obtain,
\begin{equation}
    \partial_t^{(2)}\rho\bm{u} + \bm{\nabla}\cdot{\rho}\bm{\pi}^{\rm neq} = 0,
\end{equation}
with
\begin{equation}
    {\rho}\bm{\pi}^{\rm neq} = \bm{T} - \left(\frac{1-\beta}{2\beta}\right)\left[ \rho \varsigma^2\bm{A}'\circ\left(\bm{I} \circ \bm{\nabla}\bm{u}\right) +  \bm{B}\circ\left(\bm{I} \circ {\rm diag}(\bm{\nabla}\rho)\right)\right],
\end{equation}
where $\circ$ refers to the Hadamard product. The NS viscous stress tensor is,
\begin{equation}
    \bm{T} = -\rho \nu\left(\bm{\nabla}\bm{u}+\bm{\nabla}\bm{u^\dagger} - \bm{\nabla}\cdot\bm{u}\bm{I}\right) - {\rho}\eta \bm{\nabla}\cdot\bm{u}\bm{I}, 
\end{equation}
with
\begin{align}
\nu=\pi^* \delta t\left(\frac{1}{2\beta}-\frac{1}{2}\right),\label{eq:viscosity1}
\end{align}
and
\begin{equation}\label{eq:bulk_visc_app}
    \eta = \nu \left(\frac{2+D}{D} - \frac{\partial \ln(\rho\pi^*)}{\partial\ln \rho}\right).
\end{equation}
Furthermore, $\bm{A}' = {\rm diag}(\mathcal{A}'_{xx}, \mathcal{A}'_{yy})$ with,
\begin{equation}\label{eq:error_A_2d_app}
    \mathcal{A}'_{\alpha\alpha} = -\frac{3 u_\alpha^2 + 3(\pi^* - \varsigma^2)}{\varsigma^2},
\end{equation}
and $\bm{B} = {\rm diag}(\mathcal{B}_{xx}, \mathcal{B}_{yy})$ with,
\begin{equation}\label{eq:error_B_2d_app}
    \mathcal{B}_{\alpha\alpha} = -u_\alpha^3 - 3u_\alpha(\partial_\rho (\rho \pi^*) - \varsigma^2).
\end{equation}

\section{Spectral analysis of linearized balance equations for non-ideal fluids\label{app:NS_non_ideal}}
Here we detail the linearization process of the non-ideal NS balance equations and derive spectral dispersion and dissipation of eigen-modes. To that end, we introduce the following first-order expansion around a global equilibrium state $\left(\bar{\rho},\bar{\bm{u}}\right)$,
\begin{equation}
    \rho = \bar{\rho} + \rho',\,\, \bm{u} = \bar{\bm{u}} + \bm{u}'
\end{equation}
where $\rho'$ and $\bm{u}'$ are the perturbations. For the sake of readability, we will remove the bars in the remainder of text. Introducing the expansions into the balance equation for $\rho$ and collecting linear terms,
\begin{equation}
    \partial_t \rho' + {\rho}\bm{\nabla}\cdot \bm{u}' + {\bm{u}}\cdot\bm{\nabla} \rho'= 0,
\end{equation}
while for the momentum balance equations,
\begin{multline}
    {\rho}\partial_t \bm{u}' + {\rho}{\bm{u}}\cdot\bm{\nabla}\bm{u}' + \left(\pi^*({\rho}) + {\rho}\partial_{{\rho}}\pi^*({\rho}))\right)\bm{\nabla} \rho' \\ - {\rho} \nu \bm{\nabla}\cdot\left(\bm{\nabla}\bm{u}'+{\bm{\nabla}\bm{u}'}^\dagger - \frac{2}{D}\bm{\nabla}\cdot\bm{u}'\bm{I}\right) \\- \rho\eta\bm{\nabla}\cdot \left(\bm{\nabla}\cdot\bm{u}'\bm{I}\right)= 0.
\end{multline}
Next considering the perturbations to be monochromatic plane waves of the form,
\begin{equation}
    \{\rho',\bm{u}'\} = {\{\hat{\rho},\hat{\bm{u}}\}}\exp{[\mathrm{i}(\bm{k\cdot\bm{r}}-\omega t)]},
\end{equation}
and introducing them into the linearized equations, one obtains,
\begin{widetext}
	\begin{align}
  &  {\mathrm{i}\omega \hat{\rho} - \mathrm{i}{\rho} \bm{k}\cdot\hat{\bm{u}} - \mathrm{i}\bm{u}\cdot\bm{k} \hat{\rho}= 0,}\label{eq:fourierrho}\\
  &  {\mathrm{i}\omega \hat{\bm{u}} - \mathrm{i} (\bm{k}\cdot\bm{u})\hat{\bm{u}} - \mathrm{i} \bm{k}\left(\frac{\pi^*}{{\rho}} + \partial_{{\rho}}\pi^*\right) \hat{\rho} + \nu{k}^2 \hat{\bm{u}}  + \left(\frac{D-2}{D}\nu + \eta\right) \bm{k}(\bm{k}\cdot\hat{\bm{u}}) = 0.}\label{eq:fourieru}
	  \end{align}
\end{widetext}
The system \eqref{eq:fourierrho} and \eqref{eq:fourieru} can be written as an eigen-value problem,
\begin{equation}
    \omega \Phi = {M}\Phi.
\end{equation}
Specifying for the two-dimensional case, we have $\Phi = \left[\hat{\rho}, \hat{u}_x, \hat{u}_y\right]^\dagger$ and,
\begin{widetext}
\begin{equation}
    {M} = \begin{bmatrix}
                    \bm{k}\cdot{\bm{u}} & {\rho}k_x & {\rho}k_y\\
                    \left(\frac{\pi^*}{{\rho}} + \partial_{{\rho}}\pi^*\right) k_x & \bm{k}\cdot{\bm{u}}+\mathrm{i}(\nu{k}^2 + \eta k_x^2) & \mathrm{i} \eta k_x k_y\\
                    \left(\frac{\pi^*}{{\rho}} + \partial_{{\rho}}\pi^*\right) k_y & \mathrm{i} {\eta} k_x k_y & \bm{k}\cdot{\bm{u}}+\mathrm{i}(\nu{k}^2 + \eta k_y^2)
\end{bmatrix}.
\end{equation}
\end{widetext}
Solving the eigen-value problem, $\det\left(M-\omega I\right)=0$, one obtains,
	\begin{align}
    \omega_{\rm shear} &= \bm{k}\cdot{\bm{u}} + \mathrm{i}\nu {k}^2,\\
    \omega_{\rm ac+} &= \bm{k}\cdot{\bm{u}} + {k}\sqrt{\partial_{{\rho}}({\rho}\pi^*)} + \mathrm{i}\frac{(\nu+\eta)}{2} {k}^2,\\
    \omega_{\rm ac-} &= \bm{k}\cdot{\bm{u}} - {k}\sqrt{\partial_{{\rho}}({\rho}\pi^*)} + \mathrm{i}\frac{(\nu+\eta)}{2} {k}^2,
	  \end{align}
indicating the existence of two acoustic modes and of a shear mode. 
\section*{references}
\bibliography{main}

\end{document}